\documentclass[12pt,a4paper]{article} 
\usepackage[english]{babel} 
\usepackage[utf8]{inputenc} 
\usepackage[breaklinks=true]{hyperref}
\usepackage{graphicx}
\usepackage{amssymb}
\usepackage{amsthm}
\usepackage{amsmath}
\usepackage{microtype}
\usepackage{scalefnt}
\usepackage{booktabs} 
\usepackage{multicol}         
\usepackage{here}
\usepackage{url}
\usepackage{xurl}
\usepackage{verbatim}
\usepackage{hyperref}
\usepackage[bottom=2cm, left=3cm, right=2cm, top=3cm]{geometry} %
\usepackage{multirow} 
\usepackage{breqn}
\usepackage{bm} 
\usepackage[]{natbib} 
\usepackage{float} 
\usepackage{graphicx} 
\theoremstyle{definition}

\title{The generalized Marshall-Olkin Lomax distribution with applications to AIDS and COVID-19 data}
\author{Alexsandro A. Ferreira \\
Universidade Federal de Pernambuco, Brazil. 
\\ alexsandro.ferreira.aaf@gmail.com \\
Gauss M. Cordeiro \\
Universidade Federal de Pernambuco, Brazil. 
\\ gauss@de.ufpe.br\\}
\date{}

\begin{document}
\maketitle

\begin{abstract}
\noindent The generalized Marshall-Olkin Lomax distribution is introduced, and its properties are easily obtained from those of the Lomax distribution. A regression model for censored data is proposed. The parameters are estimated through maximum likelihood, and consistency is verified by simulations. Three real datasets 
are selected to illustrate the superiority of the new models compared
to those from two well-known classes.\\

\end{abstract}

\section{Introduction}
Extended continuous distributions are essential in data analysis since it is possible to deal with different shapes in various fields. The Lomax distribution is widely applied in many areas, such as income and wealth inequality, biological sciences, lifetime and reliability, engineering, and actuarial sciences \citep{Murthy2004, corbellini2010fitting, ModifiedLomax2023}.

Recently, there has been a growing interest in exploring extensions of the Lomax distribution. These extensions are designed to enhance its flexibility and to model an even broader range of data more effectively. Some new extensions include the Marshall-Olkin 
Lomax (MOL) \citep{ghitany2007marshall}, beta Lomax (BL) \citep{rajab2013five}, transmuted exponentiated Lomax \citep{ashour2013transmuted}, Kumaraswamy Lomax (KL) \citep{shams2013kumaraswamy}, gamma Lomax \citep{Cordeiro2015gamma}, exponential Lomax \citep{el2015exponential}, Weibull Lomax (WL) \citep{tahir2015weibull}, Topp-Leone Lomax \citep{oguntunde2019topp}, and Maxwell Lomax \citep{Maxwelllomax2022} distributions, among others.

Generating new distributions from compounding ones makes sense when there is enough novelty in the generated model and new and interesting properties, and it gives superior fits to real data 
than other generators. This is the case of this article, which introduces the {\it generalized Marshall-Olkin Lomax} (GMOL) distribution, which follows from the generator proposed by \cite{chesneau2022alternative} and the Lomax distribution. The GMO generator, although recent and not yet well-known as the beta-G (B-G) \citep{Eugene2002} and Kumaraswamy-G (K-G) \citep{Cordeiro2011} classes, has significant potential. All these generators have two additional parameters. According to \citet{selim2020distributions}, the B-G and K-G generators have each produced more than 100 special distributions. Their popularity and extensive use show their robustness and flexibility in modeling diverse data sets. However, the emergence of the GMO generator represents an extremely attractive alternative in this field. This underscores the importance of exploring and understanding various generators in statistical analysis, as even those less known can offer substantial benefits in specific contexts.

The objectives of the article are to present the GMOL 
distribution, a more flexible model than the Lomax, BL, 
and KL, find a linear representation to determine its main properties, implement a regression model for censored data,  
apply the maximum likelihood method to estimate the parameters, 
and evaluate the performance of the new models through simulations and three applications to real data.

The article unfolds as follows: Section \ref{sec2} presents the GMOL distribution. Section \ref{sec3} provides a linear 
representation of its density function and properties. Section \ref{sec4} employs a regression model. Section \ref{sec5} reports some simulations. Applications to real data are discussed 
in Section \ref{sec6}, and conclusions in Section \ref{sec7}.

\section{The GMOL model}\label{sec2}

Let $\alpha \in (0,1]$, $\lambda \in [0,1]$, and $G(x)=G(x;\bm{\xi})$ be the cumulative distribution function (CDF) of any distribution with parameter vector $\bm{\xi}$. The GMO CDF, an alternative to the B-G and K-G families, can be expressed 
as (for $x \in \text{I\!R}$)
\begin{align}\label{MGOcdf}
F(x)= F(x;\alpha,\lambda,\bm{\xi})= \frac{\lambda\,G(x) + (1-\lambda)\,G(x)^2}{\alpha + (1 - \alpha)\,G(x)}\,.
\end{align}
	
The CDF and probability density function (PDF) of the Lomax distribution (with para\-me\-ters $\tau > 0$ and $\beta > 0$), 
say Lomax$(\tau,\beta)$, are give by (for $x >0$)
\begin{align}\label{Lomaxcdf}
G(x;\tau,\beta) = 1 - \left[\beta\,(\,\beta + x\,)^{-1}\right]^\tau\,,
\end{align}
and
\begin{align}\label{Lomaxpdf}
g(x;\tau,\beta) = \frac{\tau\beta^\tau}{(\beta + x)^{\tau + 1}}\,.
\end{align}

Inserting (\ref{Lomaxcdf}) in Equation (\ref{MGOcdf}) leads 
to the GMOL CDF 
\begin{align*}
F(x) = \frac{\lambda\left\{1 - \left[\beta\,(\,\beta + x\,)^{-1}\right]^\tau\right\} + (1-\lambda)\,\left\{1 - \left[\beta\,(\,\beta + x\,)^{-1}\right]^\tau\right\}^2}{\alpha + (1 - \alpha)\,\left\{1 - \left[\beta\,(\,\beta + x\,)^{-1}\right]^\tau\right\}}\,.
\end{align*}

\begin{figure}[!ht]
		\begin{center}
			\begin{minipage}[c]{0.48\linewidth}
				\centering
				(a)
				\includegraphics[width=\textwidth, height=8cm]{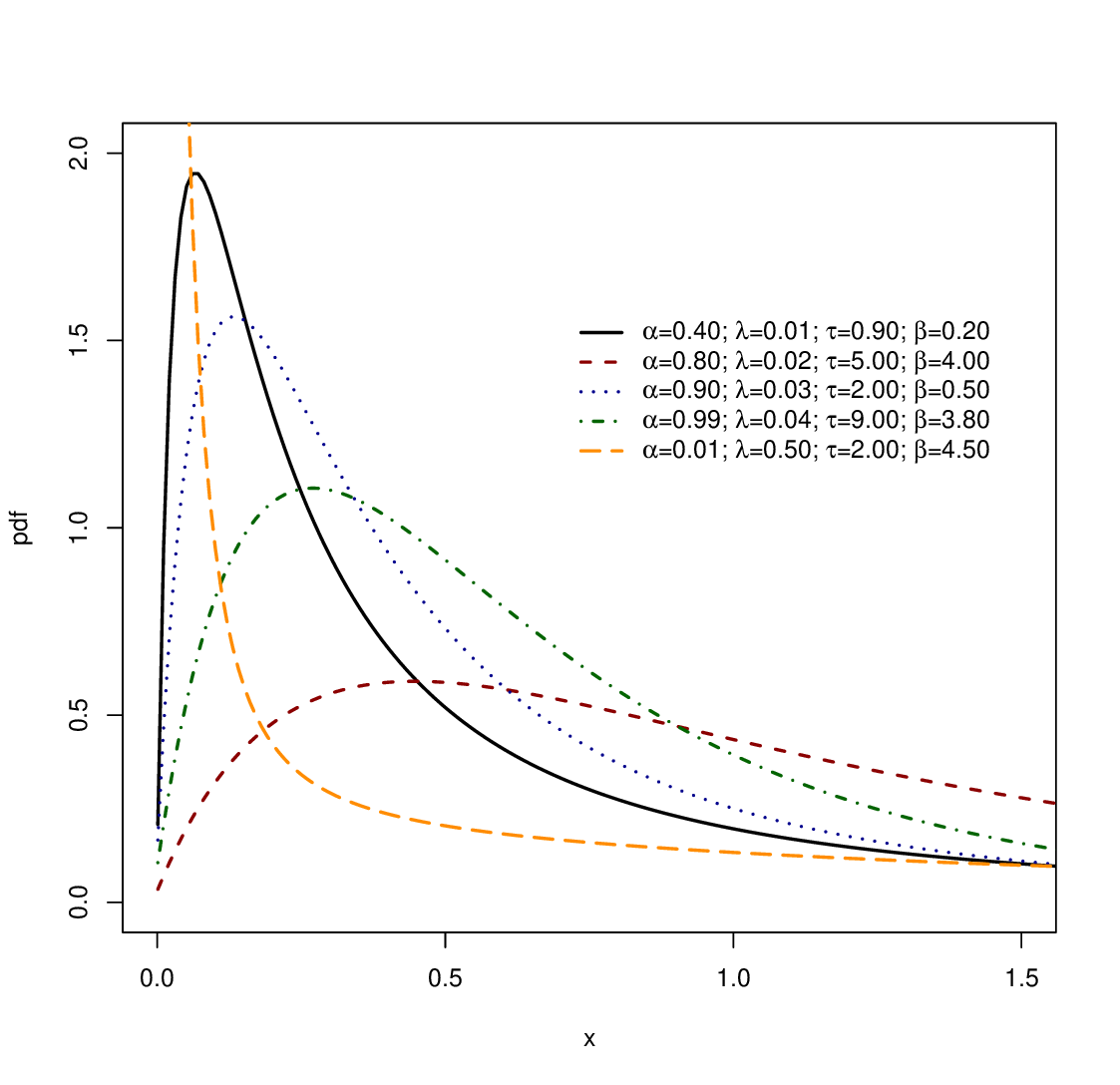}
			\end{minipage}
\hspace{.3cm}
			\begin{minipage}[c]{0.48\linewidth}
				\centering
				(b)
				\includegraphics[width=\textwidth, height=8cm]{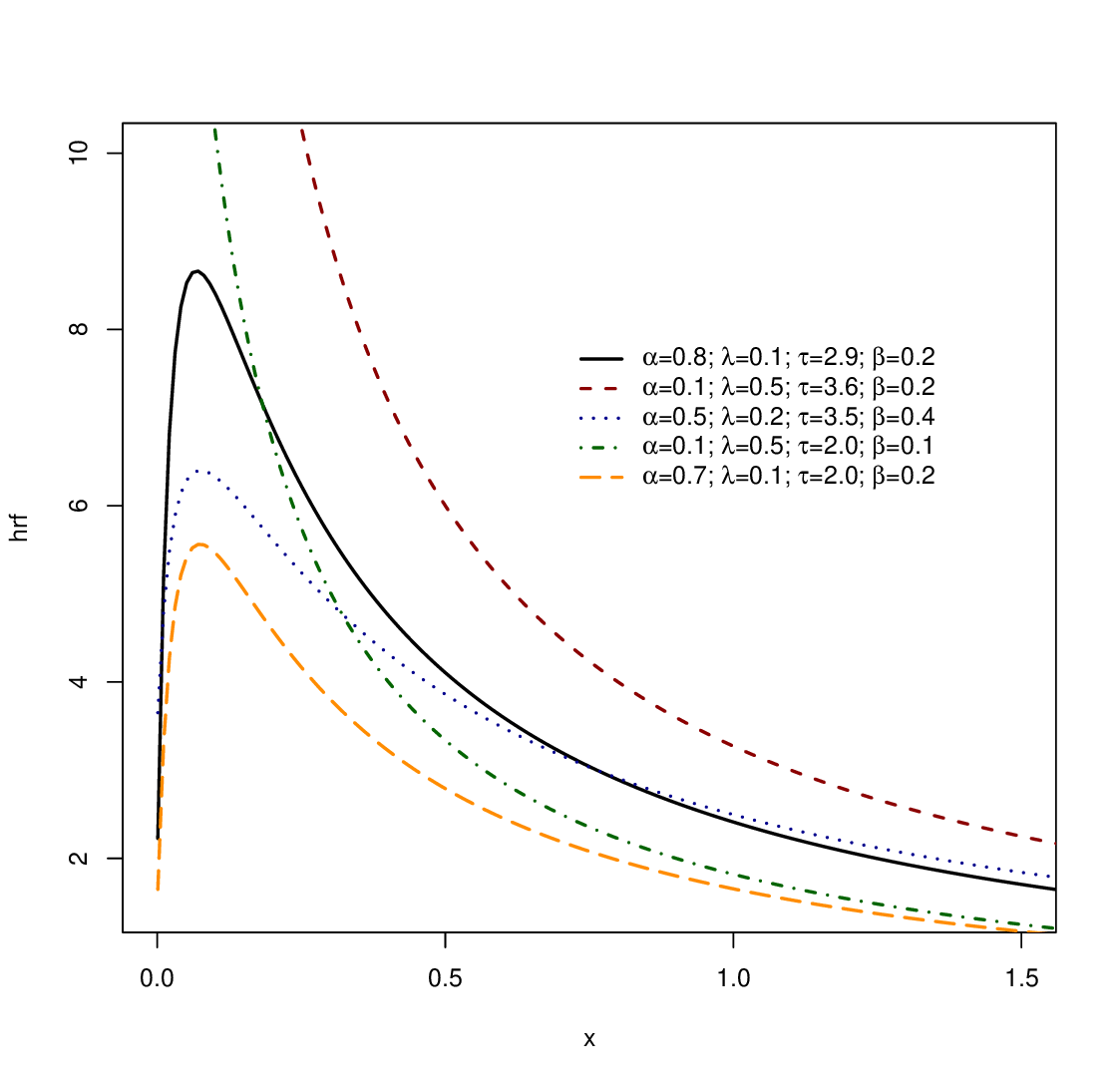}
			\end{minipage}
		\end{center}
		\caption{Density and HRF of $X$.}
	\label{F1GMOL}
	\end{figure}

The corresponding PDF and HRF can be written as
\begin{align}\label{GMOLpdf}
f(x) = \frac{\tau\beta^\tau(1-\alpha)(1-\lambda)\left\{1 - \left[\beta\,(\,\beta + x\,)^{-1}\right]^\tau\right\}^2\! + 2\alpha(1-\lambda)\left\{1 - \left[\beta\,(\,\beta + x\,)^{-1}\right]^\tau\right\} + \alpha\lambda}{(\beta + x)^{\tau + 1}\left[\alpha + (1-\alpha)\left\{1 - \left[\beta\,(\,\beta + x\,)^{-1}\right]^\tau\right\}\right]^2},
\end{align}
and
\begin{align*}
h(x) = \frac{\tau\beta^{\tau}(1-\alpha)(1-\lambda)\left[1 - \left[\beta\,(\,\beta + x\,)^{-1}\right]^\tau\right]^2 + 2\alpha(1-\lambda)\left[1 - \left[\beta\,(\,\beta + x\,)^{-1}\right]^\tau\right] + \alpha\lambda}{(\beta + x)^{\tau + 1}\left[\alpha + (1-\lambda)\left\{ 1 - \left[\beta\,(\,\beta + x\,)^{-1}\right]^\tau \right\}\right]\left[\alpha + (1-\alpha)\left\{1 - \left[\beta\,(\,\beta + x\,)^{-1}\right]^\tau\right\}\right]}
\end{align*}
Henceforth, let $X \sim \text{GMOL}(\alpha,\lambda,\tau,\beta)$ have the PDF (\ref{GMOLpdf}). The GMOL distribution accommodates some special cases. For instance, it reduces to the MOL model when $\lambda = 1$. For $|\lambda|\, \leq 1$ and $\alpha = 1$, it becomes the transmuted Lomax \citep{ashour2013transmuted}. For $\lambda = 1$ and $\alpha = 1/2$, it is identical to a special 
case of the M Lomax \citep{kumar2017new}. Finally, if $\lambda = \alpha = 1$, it corresponds to the Lomax. Figure \ref{F1GMOL} displays the density and hazard rate shapes of the GMOL model.

\section{Properties}\label{sec3}

\subsection{Linear representation}
First, the generalized binomial expansion for $\mid \! \upsilon \! \mid\, < 1$, and $q \in \text{I\!R}$, is given by 
\begin{align}\label{binomialseriesGMOL}
(1 - \upsilon)^q = \sum^\infty_{j = 0} (-1)^j \binom{q}{j} \upsilon^{\,j}\,,
\end{align}
where $\binom{q}{0}=1$ and $\binom{q}{j} = \frac{1}{j!}\prod^j_{n=1}(q -n +1)$ (for $j \geq 1$). 

The PDF associated with (\ref{MGOcdf}) follows from Theorem 1 of \cite{chesneau2022alternative} as  
\begin{align}\label{expGMOpdf}
f(x) = g_0(x) + \sum^\infty_{i = 0}\sum^{i+1}_{j = 0}\omega_{i,j}\,g_j(x)\,,
\end{align}
where $g_j(x) = (j+1)\,g(x)\,G(x)^{j}$ is the exponentiated-G density with power $(j+1)$, and
\begin{align*}
\omega_{i,j} = (-1)^j\,(\lambda-\alpha)(1-\alpha)^i\binom{i+1}{j}\,.
\end{align*}
	
Further, we can write (\ref{expGMOpdf}) as
\begin{align}\label{expGMOpdfend}
f(x) = g_0(x) + \sum^\infty_{j = 0}\rho_j\,g_j(x)\,,
\end{align}
where $\rho_j = \sum^\infty_{i = \delta_j} \omega_{i,j}$, and $\delta_0 = \delta_1=0$ and $\delta_j = j-1$ for $j \geq 2$. Thus, substituting (\ref{Lomaxcdf}) and (\ref{Lomaxpdf}) in Equation  (\ref{expGMOpdfend}), the PDF of $X$ reduces to
\begin{align}\label{expGMOLpdf}
f(x) = \frac{\tau\beta^\tau}{(\beta + x)^{\tau + 1}} + \sum^\infty_{j = 0}\rho_j\,\frac{(j+1)\tau\beta^\tau}{(\beta + x)^{\tau + 1}}\left\{1 - \left[\beta\,(\,\beta + x\,)^{-1}\right]^\tau\right\}^j\,.
\end{align}
	
Applying expansion (\ref{binomialseriesGMOL}) to the binomial term in (\ref{expGMOLpdf}), 
\begin{align}\label{expGMOLpdf1}
f(x) = \frac{\tau\beta^\tau}{(\beta + x)^{\tau + 1}} + \sum^\infty_{j = 0}\rho_j\,\frac{(j+1)\tau\beta^\tau}{(\beta + x)^{\tau + 1}}\sum^{j}_{k = 0}(-1)^k\binom{j}{k}\left[\beta(\beta + x)^{-1}\right]^{k\tau}\,.
\end{align}

By rearranging the terms in (\ref{expGMOLpdf1}) and  changing $\sum^\infty_{j = 0}\sum^j_{k=0}$ by $\sum^\infty_{k = 0} \sum^\infty_{j = k}$,  
\begin{align}\label{expGMOLpdffinal}
f(x) = g(x; \tau,\beta) + \sum^\infty_{k = 0}\,\varphi_k\,g(x;\tau^*,\beta)\,,
\end{align}
where $\tau^* = (k+1)\tau$,
\begin{align*}
\varphi_k = \frac{(-1)^k}{(k+1)}\sum^{\infty}_{j=k}(j+1)\binom{j}{k}\rho_j\,,
\end{align*}
and $g(x;\tau^*,\beta)$ is the PDF of the Lomax distribution with parameters $\tau^*$ and $\beta$. So, the GMOL properties follow from  (\ref{expGMOLpdffinal}), and those Lomax properties.

\subsection{Quantile function}
The quantile function (qf) of $X$ is \citep{chesneau2022alternative} (for $0<u<1$) 
\begin{align}\label{GMOLqf}
Q_X(u)= \beta\left(\left\{1-\left[\frac{(1 - \alpha)u - \lambda + \sqrt{[\,\lambda - (1-\alpha)u\,]^2 + 4\alpha(1-\lambda)u}}{2(1-\lambda)}\right]\right\}^{-1/\tau} - 1\right)\,.
\end{align}
 
\begin{figure}[ht!]
		\begin{center}
			\begin{minipage}[c]{0.48\linewidth}
				\centering
			\includegraphics[width=\textwidth, height=7.5cm]{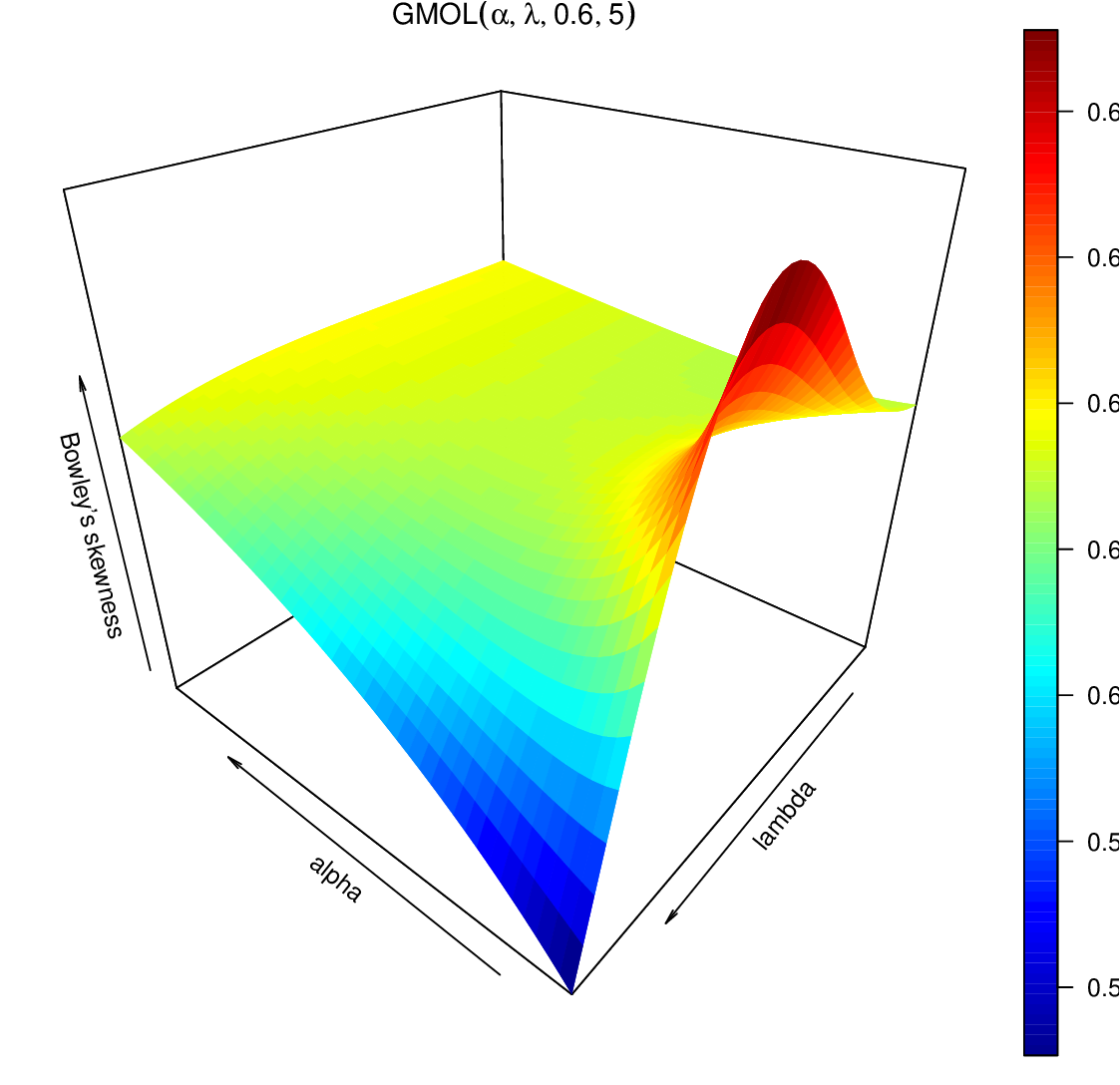}
			\end{minipage}
\hspace{.3cm}
			\begin{minipage}[c]{0.48\linewidth}
				\centering
			\includegraphics[width=\textwidth, height=7.5cm]{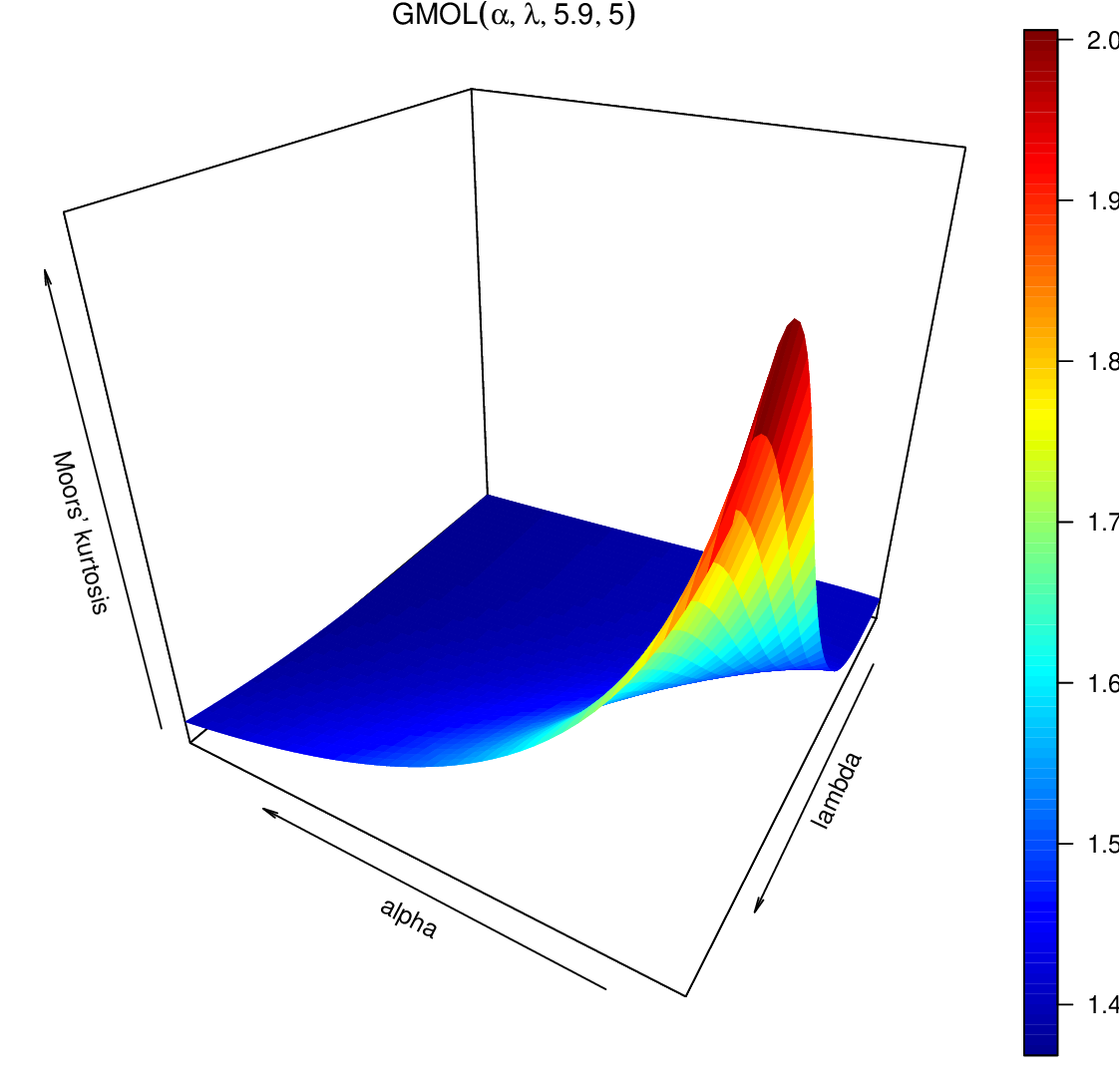}
			\end{minipage}
		\end{center}
		\caption{Skewness and kurtosis of $X$.}
	\label{skewness/kurtosis}
	\end{figure}

Plots of the Bowley skewness \citep{Kenney1962} and Moors kurtosis \citep{Moors1988} of $X$ obtained from (\ref{GMOLqf}) 
with $\tau$ and $\beta$ fixed, and varying $\alpha$ and $\lambda$, are reported in Figure \ref{skewness/kurtosis}. The skewness and kurtosis of $X$ change with the values of $\alpha$ and $\lambda$, where there is a peak (red region) indicating that these measures reach their maximum value. 

\subsection{Moments and generating function}
The $p$th moment of the Lomax model is \citep{Cordeiro2015gamma}
\begin{align}\label{Lomaxmoment}
\mu_{p,\text{Lomax}}^\prime = \frac{\beta^p\, \Gamma(\tau - p)\,\Gamma(p+1)}{\Gamma(\tau)}\,,\,\, p < \tau\,,
\end{align}
where $\Gamma(\cdot)$ is the complete gamma function. Thus, 
from Equations (\ref{expGMOLpdffinal}) and (\ref{Lomaxmoment}), 
the $p$th moment of $X$ reduces to  
\begin{align*}
\mu_p^\prime = \frac{\beta^p\, \Gamma(\tau - p)\,\Gamma(p+1)}{\Gamma(\tau)} + \sum^\infty_{k = 0}\varphi_k\,\frac{\beta^p\, \Gamma(\tau^* - p)\,\Gamma(p+1)}{\Gamma(\tau^*)}\,.
\end{align*}
	
The $p$th incomplete moment of $X$, $m_p(s) = \int^s_0 x^p f(x)\,\text{d}x$, follows from (\ref{expGMOLpdffinal}) as
\begin{align*}
m_{p}(s) = \int^s_0\frac{x^p \tau\beta^\tau}{(\beta + x)^{\tau + 1}}\text{d}x + \sum^\infty_{k = 0}\varphi_k\, \int^s_0\frac{x^p \tau^* \beta^{\tau^*}}{(\beta + x)^{\tau^* + 1}}\text{d}x\,.
\end{align*}
	
Next, by making a change of variable and using the upper incomplete beta function $B_z(a,b) = \int^{1}_z t^{a-1} (1-t)^{b-1}$, the $p$th incomplete moment of $X$ can be expressed as 
\begin{align*}
m_p(s) = \tau \beta^p B_{\beta/(\beta+s)}(\tau - p, p+1) + \sum^\infty_{k = 0}\varphi_k\, \tau^* \beta^p B_{\beta/(\beta+s)}(\tau^* - p, p+1)\,.
\end{align*}

\begin{figure}[ht!]
		\begin{center}
			\begin{minipage}[c]{0.48\linewidth}
				\centering
			\includegraphics[width=\textwidth, height=8cm]{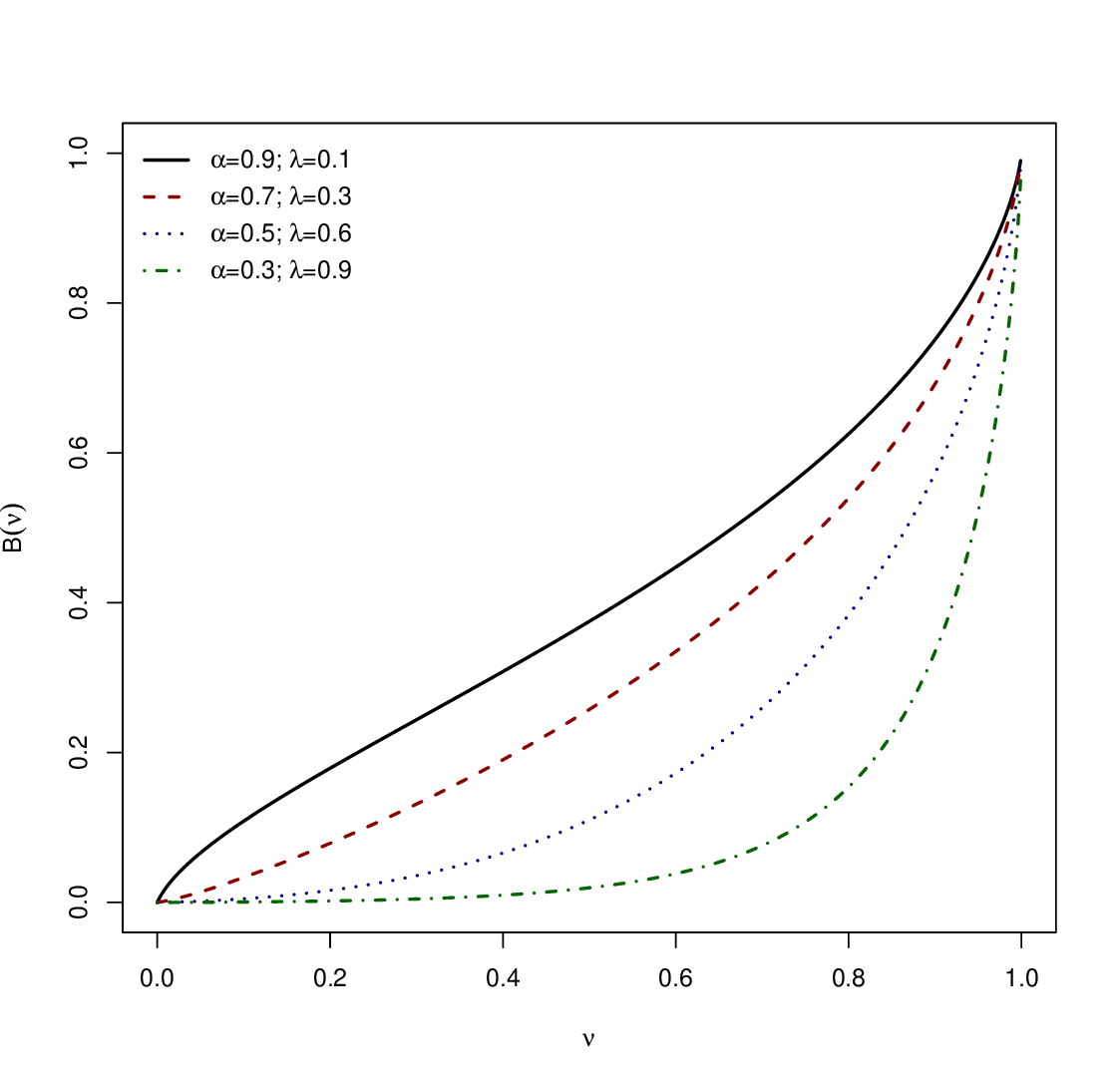}
			\end{minipage}
\hspace{.3cm}
			\begin{minipage}[c]{0.48\linewidth}
				\centering
			\includegraphics[width=\textwidth, height=8cm]{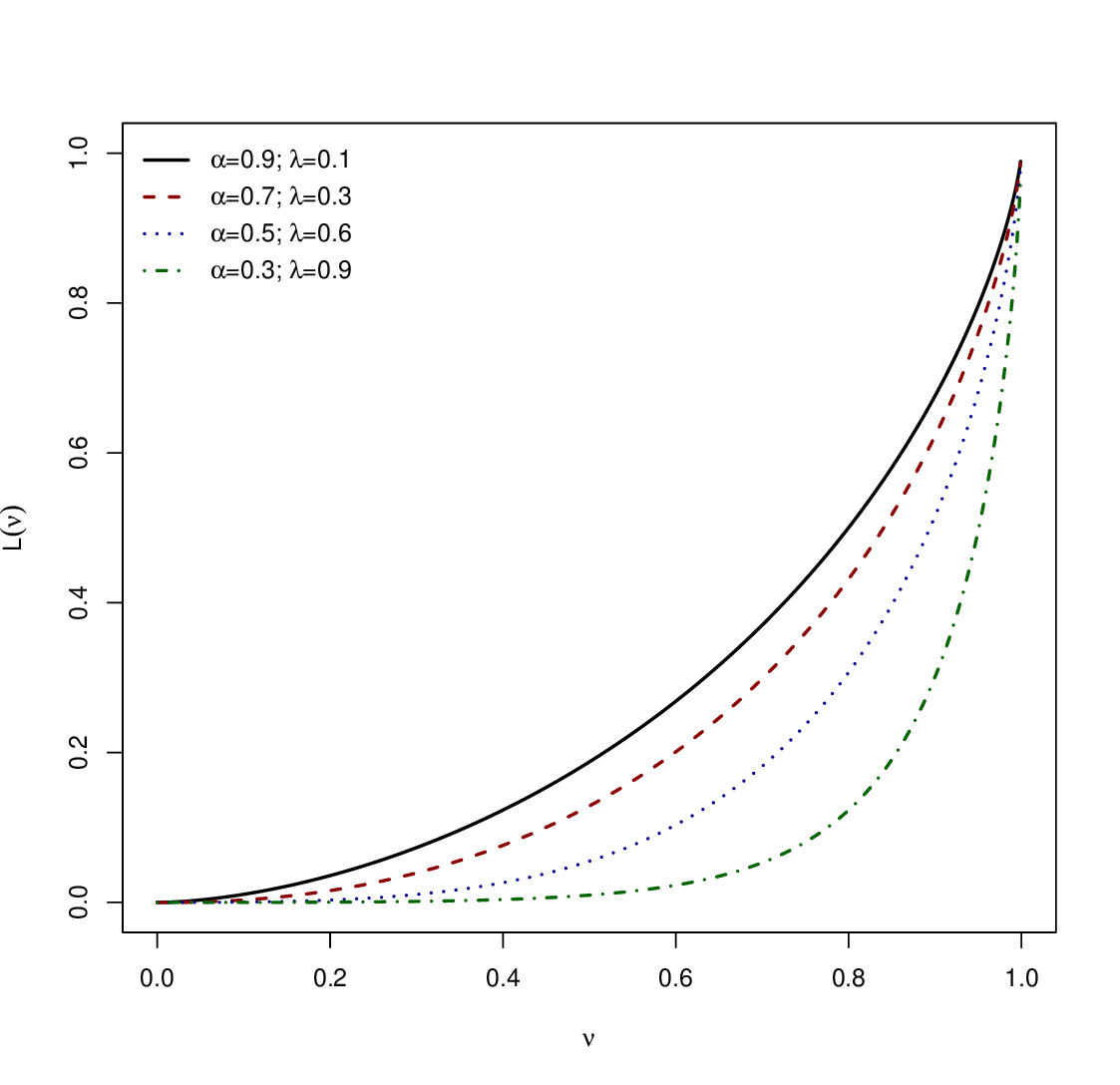}
			\end{minipage}
		\end{center}
		\caption{Bonferroni and Lorenz curves of $X$.}
	\label{bonferroni/lorenzGMOL}
\end{figure}

The first incomplete moment of $X$ provides the plots of the well-known Bonferroni and Lorenz curves. These curves versus the probability $\nu$ (varying $\alpha$ and $\lambda$) for $\tau = 6.0$ and $\beta = 2.0$, are given in Figure \ref{bonferroni/lorenzGMOL}.

Further, two analytical expressions for the generating function (gf) of $X$ are derived. First, from Equation (\ref{expGMOLpdffinal}), 
\begin{align}\label{GMOLgf}
M_X(t)=\text{I\!E}[\exp(t X)]= M_{\tau,\beta}(t) + \sum^\infty_{k = 0}\varphi_k\, M_{\tau^*,\beta}(t)\,,
\end{align}
where $M_{\tau,\beta}(t)$ and $M_{\tau^*,\beta}(t)$ are Lomax 
gfs with these parameters, namely
\begin{align*}
M_X(t) = \int^\infty_0\frac{\text{e}^{tx}\tau\beta^\tau}{(\beta + x)^{\tau + 1}}\text{d}x + \sum^\infty_{k = 0}\varphi_k\, \int^\infty_0\frac{\text{e}^{tx}\tau^* \beta^{\tau^*}}{(\beta + x)^{\tau^* + 1}}\text{d}x\,.
\end{align*}
		
For the first integral above, following the findings of \cite{lemonte2013extended}, one has
\begin{align*}
M_{\tau,\beta}(t) = \int^\infty_0\frac{\text{e}^{tx}\tau\beta^\tau}{(\beta + x)^{\tau + 1}}\text{d}x = \tau \!\int^\infty_0 \text{e}^{t\beta y}\,(1 + y)^{-(\tau + 1)}\,\text{d}y\,.
\end{align*}
	
Equations (9.210.1), (9.210.2), and (9.211.4) provided by \cite{gradshteyn2007} are employed to express the Kummer function $\Psi(a, \gamma ; z)$ as follows
\begin{align*}
\Psi(a,\gamma;z) = & \frac{1}{\Gamma(a)}\int^\infty_0 \text{e}^{-z t}\, t^{a - 1}(1+t)^{\gamma - (a+1)} \text{d}t\nonumber\\
 = & \frac{\Gamma(1-\gamma)}{\Gamma(1 + a - \gamma)}\,{}_1F_1(a,\gamma;z) + \frac{\Gamma(1-\gamma)}{\Gamma(a)}\,{}_1F_1(1 + a - \gamma,2-\gamma;z)\,,
\end{align*}
where 
\begin{align*}
{}_1F_1(a,\gamma;z) = \frac{\Gamma(\gamma)}{\Gamma(a)}\sum^\infty_{j = 0} \frac{\Gamma(a+j)z^j}{\Gamma(\gamma + j)j!}
\end{align*}
is the confluent hypergeometric function.

Then, for $t <0$, 
\begin{align}\label{Lomaxgfkummer}
M_{\tau,\beta}(t) = \Psi(1, 1 - \tau; -\beta t) = \tau^{-1}{}_1F_1(1, 1 - \tau; -\beta t) + \text{e}^{-\beta t}\,.
\end{align}

Combining (\ref{Lomaxgfkummer}) and (\ref{GMOLgf}), the gf of $X$ is
\begin{align*}
M_X(t) = \frac{1}{\tau}\,\,{}_1F_1(1, 1 - \tau; -\beta t) + \text{e}^{-\beta t} + \sum^\infty_{k = 0}\varphi_k\,\frac{1}{\tau^*}\,\,{}_1F_1(1, 1 - \tau^*; -\beta t) + \text{e}^{-\beta t}\,.
\end{align*}
	
The other representation follows from \cite{lemonte2013extended} 
by taking $u = \beta/(\beta + x)$ in (\ref{Lomaxpdf}) and solving the resulting integral (for $t < 0$) 
\begin{align}\label{lemontesecondgf}
M_{\tau,\beta}(t) = -\tau \exp(-\beta t)(-\beta t)^\tau \left[\frac{\pi(\pi\tau)}{\Gamma(\tau + 1)} + \Gamma(-\tau) - \Gamma(-\tau, -\beta t)\right]\,,\,\,t < 0\,,
\end{align}	
where $\Gamma(\cdot,\cdot)$ is the upper incomplete gamma function. Combining Equations (\ref{GMOLgf}) and (\ref{lemontesecondgf}) leads to the gf of $X$.

\subsection{Parameter estimation}

Let $x_1,\cdots,x_n$ be independent and identically 
distributed (iid) observations from $X$. The log-likelihood function for $\bm{\theta} = (\alpha,\lambda,\tau,\beta)^\top$ is
\begin{align}\label{GMOLMLE}
\ell(\bm{\theta}) =&\,\, n\log\left(\tau\beta^\tau\right) + \sum^n_{i=1}\log\Bigg\{(1-\alpha)(1-\lambda)\left\{1 - \left[\beta\,(\,\beta + x_i\,)^{-1}\right]^\tau\right\}^2\nonumber\\
& + 2\alpha(1-\lambda)\left\{1 - \left[\beta\,(\,\beta + x_i\,)^{-1}\right]^\tau\right\} + \alpha\lambda\Bigg\} - (\tau + 1)\sum^n_{i=1}\log(\beta + x_i)\nonumber\\
& - 2\sum^n_{i=1}\log\left[\alpha + (1-\alpha)\left\{1 - \left[\beta\,(\,\beta + x_i\,)^{-1}\right]^\tau\right\}\right].
\end{align}
	
The maximum likelihood estimate (MLE) of $\bm{\theta}$ can be found by maximizing (\ref{GMOLMLE}) using statistical programs. For instance, the AdequacyModel package \citep{Marinho2019} in \texttt{R} simplifies this process by providing a range of optimization methods, including Broyden-Fletcher-Goldfarb-Shannon (BFGS), Nelder-Mead, and Simulated Annealing (SANN). To use this package, one only needs to supply the respective PDFs and CDFs of the models.
	
\section{The GMOL regression model}\label{sec4}
In various real-world situations, the response variable can be affected by explanatory variables in both linear and nonlinear ways. The effects of these explanatory variables on the response variable can be evaluated through the parameters of position, scale, and shape. Recent research has introduced new regressions, such as those proposed by \cite{prataviera2018new}, \cite{PRATAVIERA2021148}, \cite{biazatti2022weibull}, and \cite{cordeirorayleigh2023}. Following a similar approach, the GMOL regression model for censored samples is defined by the systematic components $\beta_i = \exp(\bm{v}_i^\top \bm{\eta}_1)$ and $\tau_i = \exp(\bm{v}_{i}^\top \bm{\eta}_2)$ for $i = 1,\ldots,n$. Here, $\bm{v}_i^\top = (v_{i1},\cdots,v_{ir})^\top$ represents the vector of explanatory variables, while $\bm{\eta}_1 = (\eta_{11},\cdots,\eta_{1r})^\top$ and $\bm{\eta}_2 = (\eta_{21},\cdots,\eta_{2r})^\top$ denote the vectors of unknown parameters. Consequently, the survival function of $X_i | \bm{v}_i$ is 
\begin{align*}
S(x|\bm{v}_i) = 1 - \frac{\lambda\left\{1 - \left[\beta_i\,(\,\beta_i + x\,)^{-1}\right]^{\tau_i}\right\} + (1-\lambda)\,\left\{1 - \left[\beta_i\,(\,\beta_i + x\,)^{-1}\right]^{\tau_i}\right\}^2}{\alpha + (1 - \alpha)\,\left\{1 - \left[\beta_i\,(\,\beta_i + x\,)^{-1}\right]^{\tau_i}\right\}}\,.
\end{align*}
	
Consider iid observations $(x_1,\bm{v}_1),\cdots,(x_n,\bm{v}_n)$, where $x_i = \min(X_i, C_i)$, $X_i$ is the lifetime, and $C_i$ is the non-informative censoring time (assuming independence). For right-censored data, the log-likelihood function for $\bm{\zeta} = (\alpha, \lambda, \bm{\eta}_1^\top, \bm{\eta}_2^\top)^\top$ is
\begin{align}\label{GMOLmleregression}
\ell(\bm{\zeta}) =&\,\, d\log\left(\tau_i\,\beta_i^{\,\tau_i}\right) + \sum_{i\in F}\log\Bigg\{(1-\alpha)(1-\lambda)\left\{1 - \left[\beta_i\,(\,\beta_i + x_i\,)^{-1}\right]^{\tau_i}\right\}^2\nonumber\\
& + 2\alpha(1-\lambda)\left\{1 - \left[\beta_i\,(\,\beta_i + x_i\,)^{-1}\right]^{\tau_i}\right\} + \alpha\lambda\Bigg\} - (\tau_i + 1)\sum_{i\in F}\log(\beta_i + x_i) \nonumber\\
& - 2\sum_{i\in F}\log\left[\alpha + (1-\alpha)\left\{1 - \left[\beta_i\,(\,\beta_i + x\,)^{-1}\right]^{\tau_i}\right\}\right]\nonumber\\
& + \sum_{i \in C}\log\left[1 - \frac{\lambda\left\{1 - \left[\beta_i\,(\,\beta_i + x_i\,)^{-1}\right]^{\tau_i}\right\} + (1-\lambda)\,\left\{1 - \left[\beta_i\,(\,\beta_i + x_i\,)^{-1}\right]^{\tau_i}\right\}^2}{\alpha + (1 - \alpha)\,\left\{1 - \left[\beta_i\,(\,\beta_i + x_i\,)^{-1}\right]^{\tau_i}\right\}}\right]\,,
\end{align}
where $d$ is the number of failures, while $F$ and $C$ denote the lifetime and censoring sets, respectively. The MLE of $\bm{\zeta}$ can be found by numerically maximizing (\ref{GMOLmleregression}).

\section{Simulations}\label{sec5}
The MLEs of the parameters are calculated in three different scenarios with sample sizes ($n = 50$, $100$, $200$, and $300$) generated from (\ref{GMOLqf}). The average estimates (AEs), biases, and mean squared errors (MSEs) are found from one thousand Monte Carlo replicates. 
	
\begin{table}[ht!]
\centering
\small
\caption{Simulation findings from the new distribution.}
\label{TB1GMOL}
\begin{tabular}{lccccccccccc}
\\[-1.8ex] \toprule
& & \multicolumn{3}{c}{(0.2, 0.6, 0.5, 0.8)} & \multicolumn{3}{c}{(0.1, 0.3, 1.5, 3.0)} & \multicolumn{3}{c}{(0.5, 0.4, 9.0, 7.0)}\\
\cmidrule(r){3-5}\cmidrule(r){6-8}\cmidrule{9-11}
$n$ & $\bm{\theta}$ & AE & Bias & MSE & AE & Bias & MSE & AE & Bias & MSE \\
\midrule
50  & $\alpha$  & 0.2137 & 0.0136 & 0.0147 & 0.1469 & 0.0469 & 0.0271 & 1.3720  & 0.8720 & 22.2787\\
    & $\lambda$ & 0.5933 &-0.0066 & 0.0416 & 0.3702 & 0.0702 & 0.0569 & 0.7306  & 0.3306 & 3.8829 \\
    & $\tau$    & 0.5334 & 0.0334 & 0.0270 & 1.5495 & 0.0495 & 0.2730 & 10.1701 & 1.1701 & 48.8266\\
    & $\beta$   & 0.8302 & 0.0302 & 0.0268 & 3.1355 & 0.1355 & 2.1256 & 7.4160  & 0.4160 & 6.7546 \\
\\
100 & $\alpha$  & 0.2098 & 0.0098 & 0.0057 & 0.1309 & 0.0309 & 0.0137 & 0.7853  & 0.2853 & 1.2737 \\
    & $\lambda$ & 0.6110 & 0.0110 & 0.0186 & 0.3570 & 0.0570 & 0.0474 & 0.5578  & 0.1578 & 0.3929 \\
    & $\tau$    & 0.5131 & 0.0131 & 0.0112 & 1.5210 & 0.0210 & 0.0673 & 9.5443  & 0.5443 & 6.5548 \\
    & $\beta$   & 0.8230 & 0.0230 & 0.0094 & 3.1132 & 0.1132 & 0.2174 & 7.3049  & 0.3049 & 1.7086 \\
\\
200 & $\alpha$  & 0.2031 & 0.0031 & 0.0022 & 0.1219 & 0.0219 & 0.0045 & 0.6394  & 0.1394 & 0.3661 \\
    & $\lambda$ & 0.6099 & 0.0099 & 0.0067 & 0.3442 & 0.0442 & 0.0138 & 0.4931  & 0.0931 & 0.1298 \\
    & $\tau$    & 0.5076 & 0.0076 & 0.0045 & 1.5170 & 0.0170 & 0.0245 & 9.3125  & 0.3125 & 1.2109 \\
    & $\beta$   & 0.8260 & 0.0260 & 0.0039 & 3.1030 & 0.1030 & 0.0685 & 7.2276  & 0.2276 & 0.5635 \\
\\
300 & $\alpha$  & 0.2030 & 0.0030 & 0.0010 & 0.1162 & 0.0162 & 0.0026 & 0.6001  & 0.1001 & 0.1726 \\
    & $\lambda$ & 0.6128 & 0.0128 & 0.0029 & 0.3372 & 0.0372 & 0.0082 & 0.4739  & 0.0739 & 0.0790 \\
    & $\tau$    & 0.5047 & 0.0047 & 0.0023 & 1.5098 & 0.0098 & 0.0128 & 9.3029  & 0.3029 & 1.1190 \\
    & $\beta$   & 0.8255 & 0.0255 & 0.0021 & 3.1096 & 0.1096 & 0.0418 & 7.2178  & 0.2178 & 0.4575 \\
\bottomrule
\end{tabular}
\end{table}

The numbers in Table \ref{TB1GMOL} reveal that the biases 
and MSEs decay and the AEs tend to the true parameter values
if $n$ grows. This pattern confirms the consistency of the GMOL estimators.
	
One thousand Monte Carlo replicates are carried out for $n = 100, 200, 300$, and 500 to evaluate the MLEs in the regression model, whose parameters are: $\alpha=0.5$, $\lambda = 0.3$, $\eta_{10} = 0.6$, $\eta_{11} = 0.8$, $\eta_{20} = 0.2$ and $\eta_{21} = 0.4$. The censoring times $c_1,\cdots,c_n$ are generated from a uniform $(0,b)$, leading to approximately 0\%, 10\% and 30\% censoring.

The simulations follow the steps (for $i = 1,\ldots,n$):
\begin{enumerate}
\item Generate $v_{i1} \sim \text{Uniform}\,(0,1)$, and obtain $\beta_i =\exp(\eta_{10} + \eta_{11} v_{i1})$, and $\tau_i = \exp(\eta_{20} + \eta_{21} v_{i1})$.
\item  Determine the lifetimes $x_i^*$ from the GMOL$(\alpha, \lambda,\tau_i,\beta_i)$ model using Equation (\ref{GMOLqf}).
\item Generate $c_i\sim\text{Uniform}\,(0,b)$, and set $x_i = \min(x^*_i,c_i)$. If $x^*_i \leq c_i$ the censoring indicator $\delta_i = 1$. Otherwise, $\delta_i = 0$.
\end{enumerate}

\begin{table}[ht!]
\centering
\scalefont{0.9}
\caption{Simulation findings from the new regression model.}
\begin{tabular}{lccccccccccc}
\\[-1.8ex]\toprule
& & \multicolumn{3}{c}{$0\%$} & \multicolumn{3}{c}{$10\%$} & \multicolumn{3}{c}{$30\%$}\\
\cmidrule(r){3-5}\cmidrule(r){6-8}\cmidrule(r){9-11}
$n$ & $\bm{\zeta}$ & AE & Bias & MSE & AE & Bias & MSE & AE & Bias & MSE \\
\midrule
100 & $\alpha$     & 0.9329 & 0.4328 & 7.7012 & 0.9751 & 0.4751 & 9.4838 & 1.0863 & 0.5863 & 21.1799\\
    & $\lambda$    & 0.6504 & 0.3504 & 5.9446 & 0.6860 & 0.3860 & 7.3036 & 0.7851 & 0.4850 & 17.1892\\
    & $\eta_{10}$  & 0.3656 &-0.2344 & 2.3437 & 0.3609 &-0.2391 & 3.1236 & 0.3536 &-0.2463 & 4.1741 \\
    & $\eta_{11}$  & 0.6945 &-0.1055 & 2.1605 & 0.6262 &-0.1738 & 4.4126 & 0.5644 &-0.2355 & 6.8824 \\
    & $\eta_{20}$  & 0.1738 &-0.0262 & 1.3950 & 0.2038 & 0.0037 & 3.0267 & 0.1728 &-0.0272 & 4.0611 \\
    & $\eta_{21}$  & 0.4886 & 0.0886 & 0.7776 & 0.5649 & 0.1649 & 1.8629 & 0.5409 & 0.1409 & 1.2378 \\
\\
200 & $\alpha$     & 0.6389 & 0.1389 & 3.2173 & 0.6766 & 0.1765 & 4.2393 & 0.8337 & 0.3337 & 5.0693 \\ 
    & $\lambda$    & 0.4185 & 0.1185 & 2.7730 & 0.4502 & 0.1501 & 3.6030 & 0.5968 & 0.2968 & 3.3692 \\
    & $\eta_{10}$  & 0.4976 &-0.1024 & 2.6778 & 0.4783 &-0.1216 & 1.6837 & 0.3388 &-0.2612 & 1.0847 \\
    & $\eta_{11}$  & 0.7763 &-0.0237 & 0.3780 & 0.7739 &-0.0261 & 0.1149 & 0.7101 &-0.0898 & 6.5198 \\
    & $\eta_{20}$  & 0.2286 & 0.0286 & 0.0650 & 0.2120 & 0.0120 & 1.4520 & 0.2079 & 0.0079 & 2.0934 \\
    & $\eta_{21}$  & 0.4674 & 0.0674 & 0.2637 & 0.4963 & 0.0963 & 0.0988 & 0.4988 & 0.0988 & 1.0179 \\
\\
300 & $\alpha$     & 0.5813 & 0.0813 & 0.5996 & 0.5808 & 0.0808 & 0.5558 & 0.6237 & 0.1237 & 1.7251\\
    & $\lambda$    & 0.3658 & 0.0658 & 0.4516 & 0.3639 & 0.0639 & 0.4036 & 0.3639 & 0.1036 & 1.3717\\
    & $\eta_{10}$  & 0.5443 &-0.0557 & 0.4255 & 0.5503 &-0.0496 & 0.5116 & 0.5504 &-0.0883 & 1.5804\\
    & $\eta_{11}$  & 0.7868 &-0.0131 & 0.2021 & 0.7983 &-0.0017 & 0.0371 & 0.7980 &-0.0019 & 0.3882\\
    & $\eta_{20}$  & 0.2269 & 0.0269 & 0.0293 & 0.2197 & 0.0197 & 0.0427 & 0.2143 & 0.0143 & 0.0769\\
    & $\eta_{21}$  & 0.4486 & 0.0486 & 0.0522 & 0.4513 & 0.0513 & 0.0749 & 0.4574 & 0.0574 & 0.2823\\
\\
500 & $\alpha$     & 0.5305 & 0.0305 & 0.0479 & 0.5309 & 0.0309 & 0.0548 & 0.5360 & 0.0360 & 0.0699\\
    & $\lambda$    & 0.3209 & 0.0209 & 0.0283 & 0.3218 & 0.0218 & 0.0347 & 0.3249 & 0.0249 & 0.0430\\
    & $\eta_{10}$  & 0.5850 &-0.0149 & 0.0311 & 0.5851 &-0.0148 & 0.0302 & 0.5855 &-0.0144 & 0.0399\\
    & $\eta_{11}$  & 0.8058 & 0.0058 & 0.0067 & 0.8082 & 0.0082 & 0.0076 & 0.8058 & 0.0058 & 0.0111\\
    & $\eta_{20}$  & 0.2212 & 0.0212 & 0.0061 & 0.2167 & 0.0167 & 0.0100 & 0.2185 & 0.0185 & 0.0098\\
    & $\eta_{21}$  & 0.4296 & 0.0296 & 0.0080 & 0.4311 & 0.0311 & 0.0102 & 0.4353 & 0.0353 & 0.0132\\
\bottomrule \label{TB2GMOL}
\end{tabular}
\end{table}

The numbers in Table \ref{TB2GMOL} show that increasing the percentage of censoring, the biases and MSEs grow, negatively affecting their accuracy. However, larger samples mitigate these effects, even in the presence of a high censoring percentage. In general, the AEs get closer to the real parameters, and the biases and MSEs converge to zero when $n$ increases. This trend indicates the consistency of the GMOL regression estimators.
	
	Equations (\ref{GMOLMLE}) and (\ref{GMOLmleregression}) are optimized using the Nelder-Mead numerical method implemented in the \texttt{optim} function in \texttt{R} with numerical derivatives. The initial values of the parameters for both optimizations are set to their actual values. These simulation scripts are carried out without the use of any package in \texttt{R}.
	
\section{Applications}\label{sec6}
The new models are applied to three real data sets. The Cramér-von Mises $(W^*)$ and Anderson-Darling $(A^*)$ statistics \citep{chen1995}, Kolmogorov-Smirnov (KS) statistic with 
its $p$-value, as well as those well-known defined by the acronyms, are used to compare the fits of the models.  
The statistical analyses addressed in Sections \ref{sub6.1} and \ref{sub6.2} are carried out using the AdequacyModel package in \texttt{R} with the BFGS method. One of the main advantages of this package is its ability to calculate the MLEs, their standard errors (SEs), and various adequacy measures without the need to define the log-likelihood function. All that is required is to provide the PDF and CDF of the distribution that fits the data set. In Section \ref{sub6.3}, an \texttt{R} script is executed to obtain the results. This script requires the GenSA package to determine the initial values and employs the SANN numerical method via the \texttt{optim} function, whose output provides the MLEs, SEs, and adequacy measures for the chosen regressions fitted to the data.
	
\subsection{AIDS data}\label{sub6.1}
The first data set refers to the time (in days) until the occurrence of opportunistic diseases in a cohort of 695 HIV-positive patients treated at the Hospital Universitário Gaffrée e Guinle (UNIRIO) in Rio de Janeiro, Brazil, between 1995 and 2002. These data were taken from the link \url{http://sobrevida.fiocruz.br/aidsoport.html}. The average time until the occurrence of opportunistic diseases 
is 520.062 days, accompanied by a substantial standard deviation of 608.729. The skewness (1.579) and kurtosis (5.135) indicate that the data are right-skewed and leptokurtic.

The GMOL distribution is compared with well-established distributions such as BL, KL, WL, MOL, and Lomax. The MLEs and their SEs (in parentheses) for the fitted 
models are shown in Table \ref{TB3GMOL}, which reveals that all distributions, except the BL and WL distributions, provide accurate estimates.
	
\begin{figure}[ht!]
		\begin{center}
			\begin{minipage}[c]{0.48\linewidth}
				\centering
			\includegraphics[width=\textwidth, height=8cm]{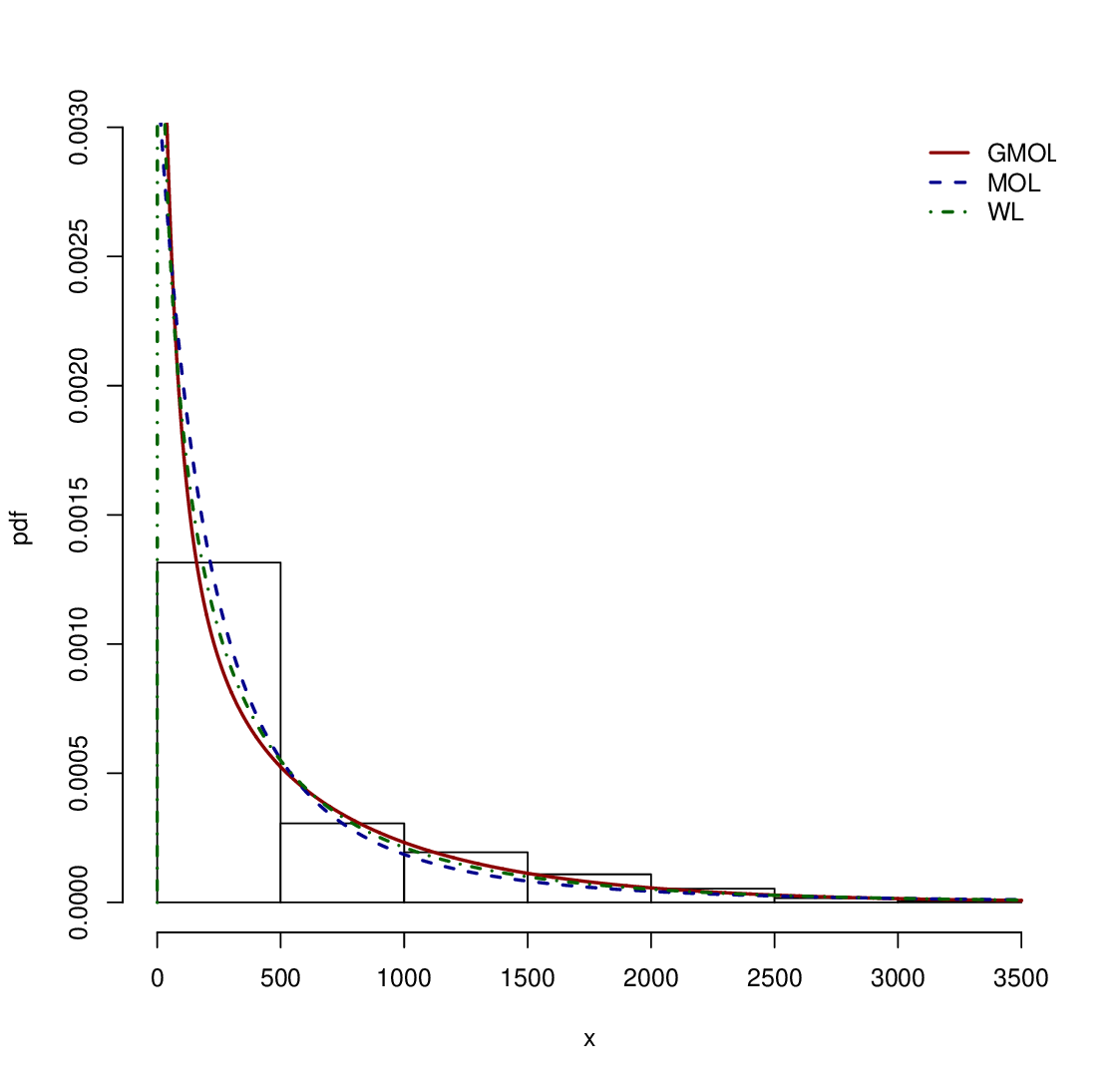}
			\end{minipage}
\hspace{.3cm}
			\begin{minipage}[c]{0.48\linewidth}
				\centering
			\includegraphics[width=\textwidth, height=8cm]{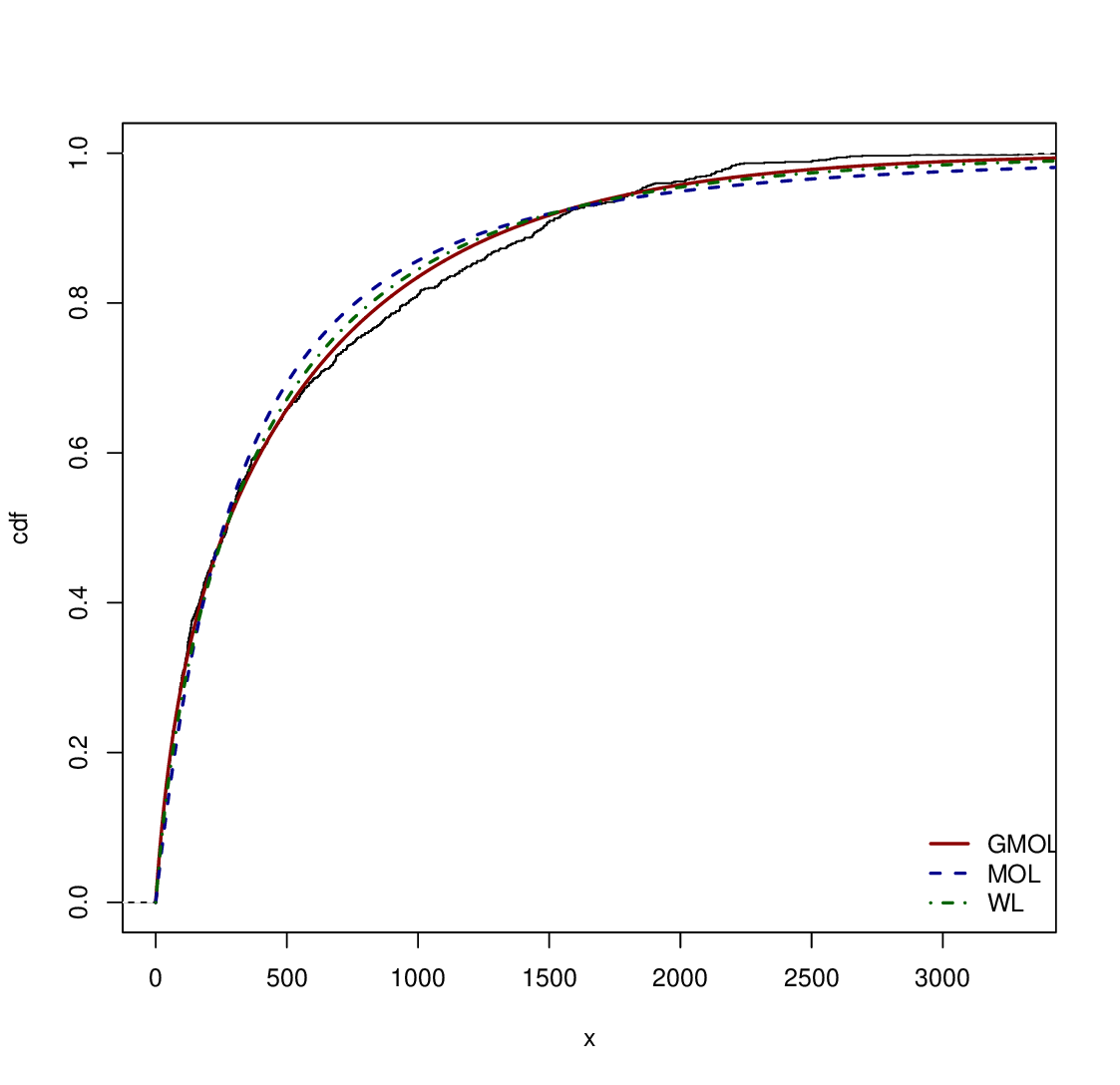}
			\end{minipage}
		\end{center}
		\caption{Estimated plots for AIDS data.}
	\label{F5GMOL}
\end{figure}

Table \ref{TB4GMOL} reveals that the GMOL distribution has 
the lowest values of the adequacy measures, thus proving a good 
fit of the model to the data when compared to the alternative ones. 
In addition, the results of the generalized likelihood ratio (GLR) tests \citep{Vuong1989} confirm the superiority of the GMOL model compared to the BL (GLR = 53.268), KL (GLR = 30.013), WL (GLR = 50.454), MOL (GLR = 35.174), and Lomax (GLR = 40.895) models at a 5\% significance level. The plots in Figure \ref{F5GMOL} also support the numerical findings.
	
\begin{table}[ht!]
\centering 
\caption{Findings for AIDS data.} 
\label{TB3GMOL} 
\begin{tabular}{@{\extracolsep{5pt}} lcccc} 
\\[-1.8ex] \toprule
\multicolumn{1}{c}{Model} & \multicolumn{4}{c}{MLEs (SEs)} \\
\midrule
GMOL & 0.137   & 0.493   & 48.753  & 34746.560 \\ 
     & (0.021) & (0.054) & (3.566) & (1861.355)\\ 
     \midrule
BL   & 0.855   & 5.920   & 0.425   & 1044.264 \\ 
     & (0.038) & (15.363)& (1.114) & (186.740)\\ 
     \midrule
KL   & 21.162  & 49.343  & 0.158   & 0.005 \\ 
     & (1.705) & (9.962) & (0.010) & (0.001)\\ 
     \midrule
WL   & 0.019   & 1.628   & 0.428   & 1.179 \\ 
     & (0.007) & (0.381) & (0.121) & (0.824)\\
     \midrule 
MOL  & 0.754   & 2.512   & 1029.416 \\ 
     & (0.127) & (0.336) & (256.111)\\ 
     \midrule
Lomax& 1.695   & 478.588 \\ 
     & (0.136) & (55.623)\\ 
\bottomrule
\end{tabular} 
\end{table} 

\begin{table}[ht!]
\scalefont{0.87}
\centering 
\caption{Adequacy measures for AIDS data.} 
\label{TB4GMOL} 
\begin{tabular}{@{\extracolsep{5pt}} lcccccccc} 
\\[-1.8ex] \toprule
 Model & $W^*$ & $A^*$ & AIC & CAIC & BIC & HQIC & KS & $p$-value    \\ 
\midrule 
GMOL & 0.255 & 1.936 & 17173.830 & 17173.860 & 17194.170 & 17181.500 & 0.029 & 0.267\\ 
BL   & 0.780 & 5.676 & 17258.240 & 17258.270 & 17278.580 & 17265.900 & 0.048 & 0.006\\ 
KL   & 0.952 & 7.004 & 17281.970 & 17282.000 & 17302.310 & 17289.630 & 0.049 & 0.007\\ 
WL   & 0.427 & 3.165 & 17199.330 & 17199.370 & 17219.670 & 17207.000 & 0.039 & 0.051\\ 
MOL  & 0.808 & 5.933 & 17268.540 & 17268.560 & 17283.800 & 17274.290 & 0.054 & 0.001\\ 
Lomax& 1.005 & 7.291 & 17291.050 & 17291.060 & 17301.220 & 17294.890 & 0.057 & $<0.001$\\ 
\bottomrule \\[-1.8ex] 
\end{tabular} 
\end{table} 

\subsection{COVID-19 data (Paraíba)}\label{sub6.2}

The second data set consists of the lifetime (in days) of 109 individuals who died from COVID-19 in the state of Paraíba, 
Brazil, from 2020 to 2022, extracted from \url{https://opendatasus.saude.gov.br/dataset/notificacoes-de-sindrome-gripal-leve-2022}. For these data, the average lifetime is 75.752 days, and the standard deviation is 100.251. The data are right-skewed and leptokurtic with skewness (1.641) and kurtosis (4.769).
	
\begin{figure}[ht!]
		\begin{center}
			\begin{minipage}[c]{0.48\linewidth}
				\centering
			\includegraphics[width=\textwidth, height=8cm]{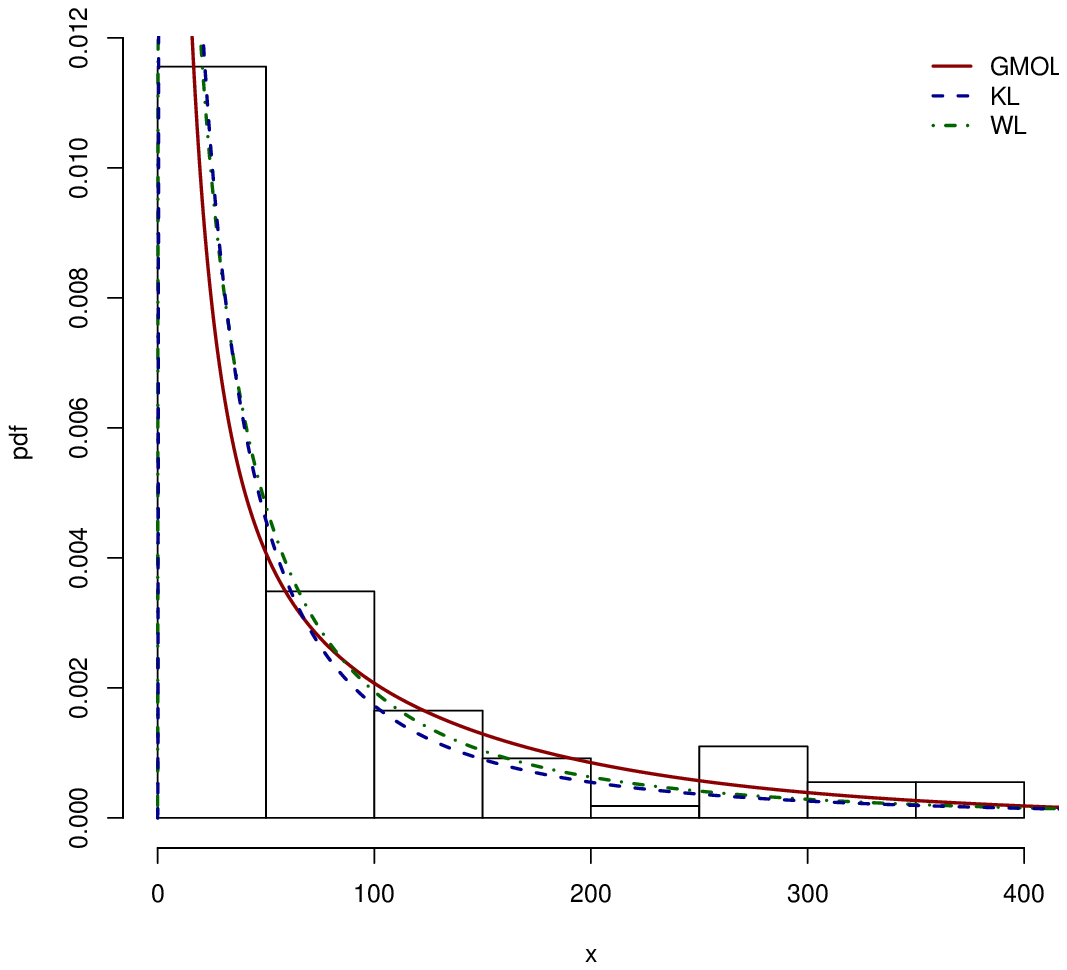}
			\end{minipage}
\hspace{.3cm}
			\begin{minipage}[c]{0.48\linewidth}
				\centering
			\includegraphics[width=\textwidth, height=8cm]{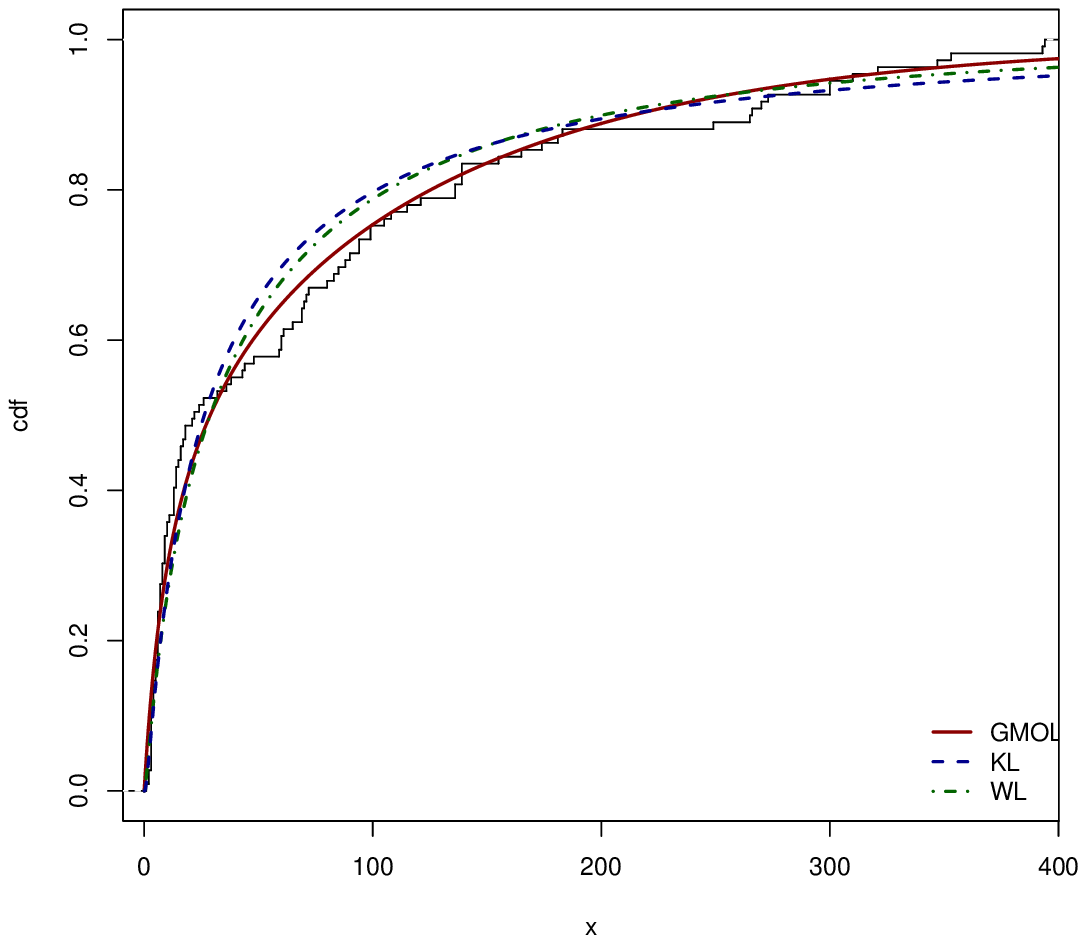}
			\end{minipage}
		\end{center}
		\caption{Estimated plots for COVID-19 data (Paraíba).}
	\label{F6GMOL}
\end{figure}

The MLEs along with their SEs for the same fitted distributions 
discussed in Section \ref{sub6.1} to the current data are provided in Table \ref{TBparaiba}. The findings in Table \ref{adequaçãoparaiba} reveal the best fit of the GMOL model to 
these data. The GLR tests also indicate that the fit 
of the GMOL model is superior to its competitors, BL (GLR = 3.609), KL (GLR = 5.395), WL (GLR = 7.995), MOL (GLR = 12.997), and Lomax (GLR = 14.908), at a 5\% significance level. These findings are visually confirmed by the plots in Figure \ref{F6GMOL}. So, both numerically and visually, the GMOL model proves to be the most flexible fitted 
distribution to these COVID-19 data.
	
\begin{table}[ht!]
\centering 
\caption{Findings for COVID-19 data (Paraíba).} 
\label{TBparaiba} 
\begin{tabular}{@{\extracolsep{5pt}} lcccc} 
\\[-1.8ex] \toprule
\multicolumn{1}{c}{Model} & \multicolumn{4}{c}{MLEs (SEs)} \\
\midrule
GMOL & 0.088  & 0.604  & 65.508 & 8674.973          \\ 
     & (0.028)  & (0.101)  & (28.619)  & (3528.939) \\ 
     \midrule
BL   & 16.852 & 12.604 & 0.136  & 0.045             \\ 
     & (13.933) & (17.752) & (0.142)  & (0.108)     \\
     \midrule 
KL   & 6.015  & 23.973 & 0.149  & 0.116             \\ 
     & (5.400)  & (40.052) & (0.117)  & (0.288)     \\ 
     \midrule
WL   & 1.701  & 2.199  & 0.132  & 0.628             \\ 
     & (8.636)  & (1.466)  & (0.245)  & (1.004)     \\ 
     \midrule
MOL  & 6.698  & 1.016  & 4.164                      \\ 
     & (0.223)  & (0.846)  & (130.195)              \\ 
     \midrule
Lomax& 1.235  & 37.244                              \\ 
     & (0.324)  & (16.420)                          \\ 
\bottomrule
\end{tabular} 
\end{table} 

\begin{table}[ht!]
\centering 
\small
\caption{Adequacy measures for COVID-19 data (Paraíba).} 
\label{adequaçãoparaiba} 
\begin{tabular}{@{\extracolsep{5pt}} lcccccccc} 
\\[-1.8ex] \toprule
 Model & $W^*$ & $A^*$ & AIC & CAIC & BIC & HQIC & KS & $p$-value    \\ 
\midrule 
GMOL & 0.259 & 1.571 & 1128.321 & 1128.706 & 1139.087 & 1132.687 & 0.098 & 0.238 \\ 
BL   & 0.320 & 1.893 & 1135.400 & 1135.785 & 1146.165 & 1139.766 & 0.130 & 0.047 \\ 
KL   & 0.328 & 1.902 & 1134.235 & 1134.620 & 1145.001 & 1138.601 & 0.115 & 0.108 \\ 
WL   & 0.326 & 1.880 & 1132.658 & 1133.043 & 1143.424 & 1137.024 & 0.102 & 0.204 \\ 
MOL  & 0.415 & 2.411 & 1141.064 & 1141.293 & 1149.138 & 1144.338 & 0.116 & 0.103 \\ 
Lomax& 0.430 & 2.503 & 1141.193 & 1141.306 & 1146.576 & 1143.376 & 0.112 & 0.127 \\ 
\bottomrule \\[-1.8ex] 
\end{tabular} 
\end{table} 

\subsection{COVID-19 data (Distrito Federal)}\label{sub6.3}

The GMOL regression model is applied to a data set composed of 485 lifetimes (in days) of COVID-19 patients in Distrito Federal, Brazil, in 2020. This data can be accessed at the link: \url{https://covid19.ssp.df.gov.br/extensions/covid19/covid19.html}, which shows an average lifetime of 8.226 days with a standard deviation of 6.219. The skewness (1.155) and kurtosis (3.992) values show that the data are right-skewed and leptokurtic.
	
The response variable $x_i$ corresponds to the time elapsed from the onset of symptoms until death from COVID-19 (failure). Approximately 78\% of the observations are censored. The variables considered are (for $i = 1,\ldots,485$): $\delta_i\!:$ censoring indicator (0 = censored, 1 = lifetime observed), $v_{i1}\!:$ age (in years), and $v_{i2}\!:$ obesity (1 = yes, 0 = no or not informed).

 In the GMOL model and its associated sub-models, namely MOL and Lomax, the explanatory variables are linked to systematic components as (for $i = 1,\ldots,485$)
\begin{align*}
\beta_i = \exp(\eta_{10} + \eta_{11}v_{i1} + \eta_{12}v_{i2})\,\, \text{and} \,\,\tau_i = \exp(\eta_{20} + \eta_{21}v_{i1} + \eta_{22}v_{i2}).
\end{align*}

According to the numbers in Table \ref{tab:criteriosGMOL}, 
the GMOL model is the most suitable of the three fitted regression models to adjust COVID-19 data in Distrito Federal. In addition, the likelihood ratio (LR) tests support the superiority of the GMOL model over the others, as indicated by the rejection of the null hypotheses (Table \ref{TBLRGMOL}).

\begin{figure}[ht!]
		\begin{center}
			\begin{minipage}[c]{0.48\linewidth}
				\centering
			\includegraphics[width=\textwidth, height=8cm]{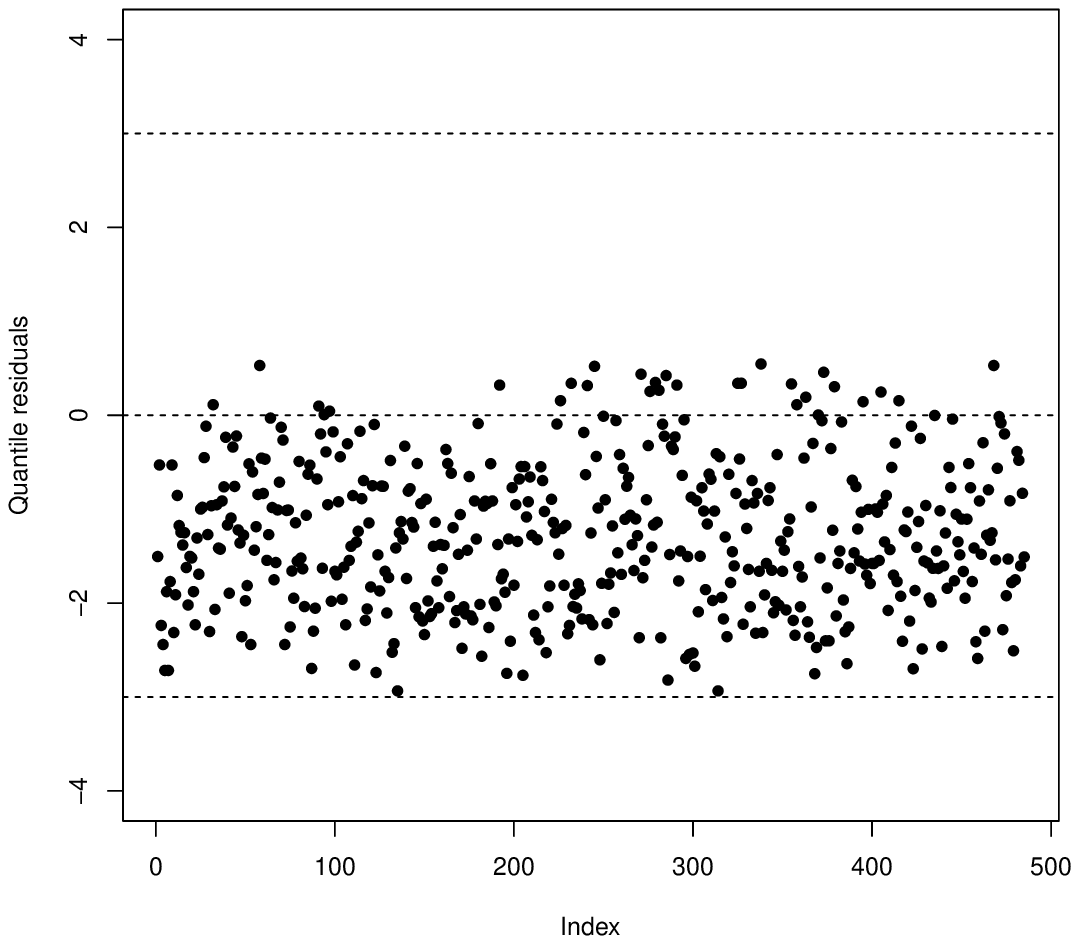}
			\end{minipage}
\hspace{.3cm}
			\begin{minipage}[c]{0.48\linewidth}
				\centering
			\includegraphics[width=\textwidth, height=8cm]{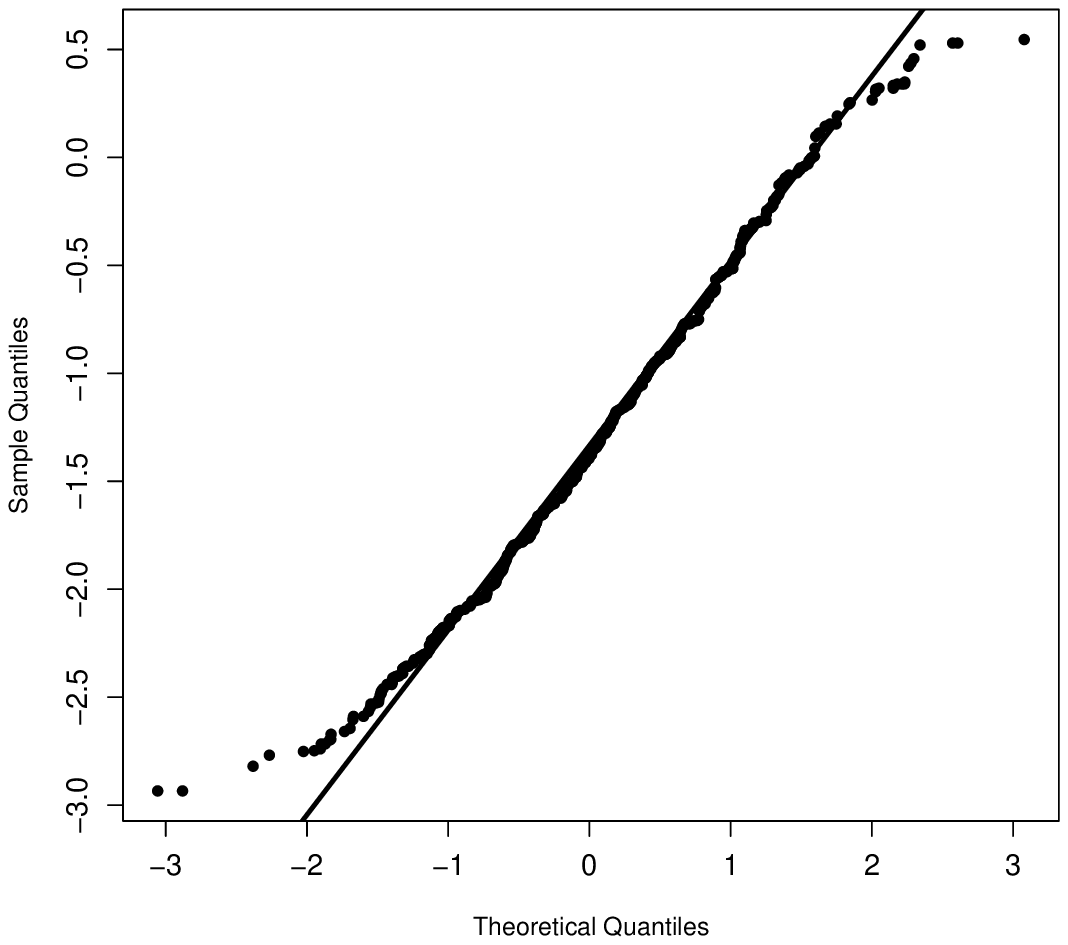}
			\end{minipage}
		\end{center}
		\caption{Index plot and normal probability plot for COVID-19 data (Distrito Federal).}
	\label{data1}
\end{figure}

The explanatory variables $v_{i1}$ ($\eta_{11}$) and $v_{i2}$ ($\eta_{12}$) are significant at the 5\% level (Table \ref{TB8GMOL}), thus indicating that age and obesity are factors that can reduce the time to failure. In addition, both are significant for the variability of survival times ($\eta_{21}$) and ($\eta_{22}$). The residual analysis for the fitted GMOL regression model is done by the quantile residuals (qrs) \citep{dunn1996}
\begin{align*}
qr_i = \Phi^{-1}\left( \frac{\hat{\lambda}\left\{1 - \left[\hat{\beta}_i\,(\,\hat{\beta}_i + x\,)^{-1}\right]^{\hat{\tau}_i}\right\} + (1-\hat{\lambda})\,\left\{1 - \left[\hat{\beta}_i\,(\,\hat{\beta}_i + x\,)^{-1}\right]^{\hat{\tau}_i}\right\}^2}{\hat{\alpha} + (1 - \hat{\alpha})\,\left\{1 - \left[\hat{\beta}_i\,(\,\hat{\beta}_i + x\,)^{-1}\right]^{\hat{\tau}_i}\right\}}\right)\,,
\end{align*}
where $\Phi(\cdot)^{-1}$ is the qf of the standard normal distribution, $\hat{\beta}_i = \exp(\bm{v}_i^\top \hat{\bm{\eta}_1})$ and $\hat{\tau}_i = \exp(\bm{v}_{i}^\top \hat{\bm{\eta}}_2)$. Figure \ref{data1} indicates a well-fitted regression model since the qrs have random behavior and approximately follow a standard normal distribution.

\begin{table}[ht!]
\centering
\caption{Adequacy measures for COVID-19 data (Distrito Federal).}
\label{tab:criteriosGMOL}
\begin{tabular}{@{\extracolsep{5pt}} lcccc} 
\\[-1.8ex]\toprule
Model & AIC & CAIC & BIC & HQIC \\
\midrule
GMOL  & 937.022 & 937.486 & 970.495 & 950.174 \\
MOL   & 942.541 & 942.920 & 971.830 & 954.048 \\
Lomax & 943.274 & 943.576 & 968.379 & 953.137 \\
\bottomrule
\end{tabular}
\end{table}

\begin{table}[ht!]
\centering
\caption{LR tests for COVID-19 data (Distrito Federal).}
\label{TBLRGMOL}
\begin{tabular}{@{\extracolsep{5pt}} ccccc} 
\\[-1.8ex]\toprule
Model & Hypotheses & LR statistic & $p$-value \\
\midrule
GMOL vs MOL   & $H_0 :\lambda=1$ vs $H_1 : H_0$ is false      & 7.5186  & 0.0061 \\
GMOL vs Lomax & $H_0 :\alpha=\lambda=1$ vs $H_1 : H_0$ is false & 10.2516 & 0.0059 \\
\bottomrule
\end{tabular}
\end{table}

\begin{table}[ht!]
\centering
\caption{Estimated results for COVID-19 data (Distrito Federal).} 
\label{TB8GMOL} 
\begin{tabular}{@{\extracolsep{5pt}} cccccccccc} 
\\[-1.8ex]\toprule
\multicolumn{1}{c}{Parameter} & \multicolumn{1}{c}{MLEs} & \multicolumn{1}{c}{SEs} & \multicolumn{1}{c}{$p$-value} \\
\midrule 
$\alpha$    & 0.967 & 0.274 & $-$      \\
$\lambda$   & 0.143 & 0.069 & $-$      \\
$\eta_{10}$ &10.248 & 2.010 & $<$0.001 \\
$\eta_{11}$ &-0.097 & 0.025 & 0.001    \\
$\eta_{12}$ &-8.661 & 3.129 & 0.005    \\
$\eta_{20}$ & 4.287 & 1.852 & 0.021    \\
$\eta_{21}$ &-0.048 & 0.021 & 0.027    \\
$\eta_{22}$ &-4.441 & 0.866 & $<$0.001 \\
\bottomrule
\end{tabular}
\end{table}

\section{Conclusions}\label{sec7}
This article presented the generalized Marshall-Olkin Lomax (GMOL) distribution, a more flexible alternative to the Lomax model. Its properties can be easily found from a linear representation in terms of Lomax densities. A regression model with two systematic components is proposed for censored data. Simulations confirm the consistency of the maximum likelihood estimators and the effectiveness of the regression model for censored data. 
The GMOL distribution outperformed the beta Lomax and 
Kumaraswamy Lomax models in real data analysis. Thus, this 
research contributed to the ongoing development of lifetime models and offered valuable tools for many diverse types of data. 

\section*{Acknowledgments}
Fundação de Amparo à Ciência e Tecnologia do Estado 
de Pernambuco (FACEPE) [IBPG-1448-1.02/20] supports this work. 

\bibliographystyle{apalike}
\bibliography{referencias}

\begin{thebibliography}{}

\bibitem[Abiodun and Ishaq, 2022]{Maxwelllomax2022}
Abiodun, A.~A. and Ishaq, A.~I. (2022).
\newblock {On Maxwell–Lomax distribution: properties and applications}.
\newblock {\em Arab Journal of Basic and Applied Sciences}, 29:221--232.

\bibitem[Alnssyan, 2023]{ModifiedLomax2023}
Alnssyan, B. (2023).
\newblock {The Modified-Lomax Distribution: Properties, Estimation Methods, and
  Application}.
\newblock {\em Symmetry}, 15:1367.

\bibitem[Ashour and Eltehiwy, 2013]{ashour2013transmuted}
Ashour, S. and Eltehiwy, M. (2013).
\newblock {Transmuted exponentiated Lomax distribution}.
\newblock {\em Australian Journal of Basic and Applied Sciences}, 7:658--667.

\bibitem[Biazatti et~al., 2022]{biazatti2022weibull}
Biazatti, E.~C., Cordeiro, G.~M., Rodrigues, G.~M., Ortega, E. M.~M., and
  de~Santana, L.~H. (2022).
\newblock A {W}eibull-beta prime distribution to model {COVID}-19 data with the
  presence of covariates and censored data.
\newblock {\em Stats}, 5(4):1159--1173.

\bibitem[Chen and Balakrishnan, 1995]{chen1995}
Chen, G. and Balakrishnan, N. (1995).
\newblock {A general purpose approximate goodness-of-fit test}.
\newblock {\em Journal of Quality Technology}, 27:154--161.

\bibitem[Chesneau et~al., 2022]{chesneau2022alternative}
Chesneau, C., Karakaya, K., Bakouch, H.~S., and Ku{\c{s}}, C. (2022).
\newblock {An alternative to the Marshall-Olkin family of distributions:
  bootstrap, regression and applications}.
\newblock {\em Communications on Applied Mathematics and Computation},
  4:1229–1257.

\bibitem[Corbellini et~al., 2010]{corbellini2010fitting}
Corbellini, A., Crosato, L., Ganugi, P., and Mazzoli, M. (2010).
\newblock {Fitting Pareto II distributions on firm size: Statistical
  methodology and economic puzzles}.
\newblock {\em Advances in Data Analysis: Theory and Applications to
  Reliability and Inference, Data Mining, Bioinformatics, Lifetime Data, and
  Neural Networks}, pages 321--328.

\bibitem[Cordeiro and de~Castro, 2011]{Cordeiro2011}
Cordeiro, G.~M. and de~Castro, M. (2011).
\newblock A new family of generalized distributions.
\newblock {\em Journal of Statistical Computation and Simulation}, 81:883--898.

\bibitem[Cordeiro et~al., 2015]{Cordeiro2015gamma}
Cordeiro, G.~M., Ortega, E. M.~M., and Popović, B.~V. (2015).
\newblock {The gamma-Lomax distribution}.
\newblock {\em Journal of Statistical Computation and Simulation}, 85:305--319.

\bibitem[Cordeiro et~al., 2023]{cordeirorayleigh2023}
Cordeiro, G.~M., Rodrigues, G.~M., Ortega, E. M.~M., de~Santana, L.~H., and
  Vila, R. (2023).
\newblock An extended {R}ayleigh model: Properties, regression and covid-19
  application.
\newblock {\em Chilean Journal of Statistics}, 14(1):1--25.

\bibitem[Dunn and Smyth, 1996]{dunn1996}
Dunn, P.~K. and Smyth, G.~K. (1996).
\newblock Randomized quantile residuals.
\newblock {\em Journal of Computational and Graphical Statistics}, 5:236--244.

\bibitem[El-Bassiouny et~al., 2015]{el2015exponential}
El-Bassiouny, A., Abdo, N., and Shahen, H. (2015).
\newblock {Exponential Lomax distribution}.
\newblock {\em International Journal of Computer Applications}, 121:24--29.

\bibitem[Eugene et~al., 2002]{Eugene2002}
Eugene, N., Lee, C., and Famoye, F. (2002).
\newblock Beta-normal distribution and its applications.
\newblock {\em Communications in Statistics - Theory and Methods}, 31:497--512.

\bibitem[Ghitany et~al., 2007]{ghitany2007marshall}
Ghitany, M.~E., Al-Awadhi, F.~A., and Alkhalfan, L. (2007).
\newblock {Marshall--Olkin extended Lomax distribution and its application to
  censored data}.
\newblock {\em Communications in Statistics—Theory and Methods},
  36:1855--1866.

\bibitem[Gradshteyn and Ryzhik, 2007]{gradshteyn2007}
Gradshteyn, I.~S. and Ryzhik, I.~M. (2007).
\newblock {\em Table of integrals, series, and products}.
\newblock Elsevier/Academic Press, Amsterdam, seventh edition.
\newblock Translated from the Russian, Translation edited and with a preface by
  Alan Jeffrey and Daniel Zwillinger, With one CD-ROM (Windows, Macintosh and
  UNIX).

\bibitem[Kenney and Keeping, 1962]{Kenney1962}
Kenney, J. and Keeping, E. (1962).
\newblock {\em Moving averages}.
\newblock 3 edn. NJ: Van Nostrand.

\bibitem[Kumar et~al., 2017]{kumar2017new}
Kumar, D., Singh, U., Singh, S.~K., and Mukherjee, S. (2017).
\newblock {The new probability distribution: an aspect to a life time
  distribution}.
\newblock {\em Math. Sci. Lett}, 6:35--42.

\bibitem[Lemonte and Cordeiro, 2013]{lemonte2013extended}
Lemonte, A.~J. and Cordeiro, G.~M. (2013).
\newblock {An extended Lomax distribution}.
\newblock {\em Statistics}, 47:800--816.

\bibitem[Marinho et~al., 2019]{Marinho2019}
Marinho, P. R.~D., Silva, R.~B., Bourguignon, M., Cordeiro, G.~M., and
  Nadarajah, S. (2019).
\newblock {AdequacyModel: An R package for probability distributions and
  general purpose optimization}.
\newblock {\em PLOS ONE}, 14:1--30.

\bibitem[Moors, 1988]{Moors1988}
Moors, J. J.~A. (1988).
\newblock A quantile alternative for kurtosis.
\newblock {\em Journal of the Royal Statistical Society. Series D (The
  Statistician)}, 37:25--32.

\bibitem[Murthy et~al., 2004]{Murthy2004}
Murthy, D.~P., Xie, M., and Jiang, R. (2004).
\newblock {\em {Weibull models}}, volume 505.
\newblock John Wiley and Sons.

\bibitem[Oguntunde et~al., 2019]{oguntunde2019topp}
Oguntunde, P.~E., Khaleel, M.~A., Okagbue, H.~I., and Odetunmibi, O.~A. (2019).
\newblock {The Topp--Leone Lomax (TLLo) distribution with applications to
  airbone communication transceiver dataset}.
\newblock {\em Wireless Personal Communications}, 109:349--360.

\bibitem[Prataviera et~al., 2018]{prataviera2018new}
Prataviera, F., Ortega, E. M.~M., Cordeiro, G.~M., Pescim, R.~R., and Verssani,
  B. A.~W. (2018).
\newblock A new generalized odd log-logistic flexible {W}eibull regression
  model with applications in repairable systems.
\newblock {\em Reliability Engineering \& System Safety}, 176:13--26.

\bibitem[Prataviera et~al., 2021]{PRATAVIERA2021148}
Prataviera, F., Silva, A. M.~M., Cardoso, E. J. B.~N., Cordeiro, G.~M., and
  Ortega, E. M.~M. (2021).
\newblock A novel generalized odd log-logistic {M}axwell-based regression with
  application to microbiology.
\newblock {\em Applied Mathematical Modelling}, 93:148--164.

\bibitem[Rajab et~al., 2013]{rajab2013five}
Rajab, M., Aleem, M., Nawaz, T., and Daniyal, M. (2013).
\newblock {On five parameter beta Lomax distribution}.
\newblock {\em Journal of Statistics}, 20:102--118.

\bibitem[Selim, 2020]{selim2020distributions}
Selim, M.~A. (2020).
\newblock {The Distributions of Beta-Generated and Kumaraswamy-Generalized
  Families: A Brief Survey}.
\newblock {\em Figshare}, pages 1--19.

\bibitem[Shams, 2013]{shams2013kumaraswamy}
Shams, T.~M. (2013).
\newblock {The Kumaraswamy-generalized lomax distribution}.
\newblock {\em Middle-East Journal of Scientific Research}, 17:641--646.

\bibitem[Tahir et~al., 2015]{tahir2015weibull}
Tahir, M.~H., Cordeiro, G.~M., Mansoor, M., and Zubair, M. (2015).
\newblock {The Weibull-Lomax distribution: properties and applications}.
\newblock {\em Hacettepe Journal of Mathematics and Statistics}, 44:455--474.

\bibitem[Vuong, 1989]{Vuong1989}
Vuong, Q.~H. (1989).
\newblock Likelihood ratio tests for model selection and non-nested hypotheses.
\newblock {\em Econometrica}, 57:307--333.

\end{thebibliography}
\end{document}